\documentclass[openany,oneside,a4paper,11pt,article]{memoir}

\usepackage[utf8]{inputenc}
\usepackage[T1]{fontenc}
\usepackage{amssymb}
\usepackage[round]{natbib}

\usepackage[french,italian,german,dutch,latin,english]{babel}

\usepackage[nocritical,noend,noledgroup,series={A,B,C,D}]{reledmac}
\usepackage[shiftedpstarts]{reledpar}
\usepackage[shiftedpstarts]{reledpar}
\usepackage{url}
\usepackage{alphalph}
\usepackage[noorphans,leftmargin=3em,rightmargin=3em]{quoting}

\usepackage{indentfirst}

\isopage

\feetatbottom

\DeclareRobustCommand{\VAN}[3]{#2}

\newcommand{\german}[1]{\foreignlanguage{german}{#1}}

\newcommand{\dutch}[1]{\foreignlanguage{dutch}{#1}}
\newcommand{\french}[1]{\foreignlanguage{french}{#1}}

\newcommand{\italian}[1]{\foreignlanguage{italian}{#1}}

\numberlinefalse% turn off eledmac's line numbering

\rightnoteupfalse
\leftnoteupfalse

\newcommand{\pagenumber}[1]{\mbox{\textbar}\ledrightnote{#1}}

\setgoalfraction{0.9}

\begin{document}

\thanksmarkseries{alph}
\setlength{\thanksmarkwidth}{0em}
%\continuousmarks

\title{L.E.J.~Brouwer's `Unreliability of the logical principles'.\\ A new translation, with an introduction}

\author{Mark van Atten\thanks{SND (CNRS/Paris 4), 28 rue Serpente, 75006 Paris, France. vanattenmark@gmail.com}
\and
Göran Sundholm\thanks{Philosophical Institute, Leiden University, P.O.~Box 2315,  2300 RA Leiden, The Netherlands.\protect\par \hspace{\parindent}goran.sundholm@gmail.com}}

\maketitle

\begin{center}
\itshape Dedicated to the memory of Georg Kreisel, 1923–2015
\end{center}

\begin{abstract}
We present a new English translation of L.E.J.~Brouwer's paper
`\dutch{De onbetrouwbaarheid der logische principes}'
(The unreliability of the logical principles)
of 1908,
together with a philosophical and historical introduction.
In this paper Brouwer for the first time
objected to the idea that the Principle of the Excluded Middle
is valid.
We discuss the circumstances under which the manuscript
was submitted and accepted,
Brouwer’s ideas on the principle of the excluded middle,
its consistency and partial validity,
and his argument against the possibility of
absolutely undecidable propositions.
We note that principled objections to 
the general excluded middle 
similar to Brouwer's
had been advanced in print
by Jules Molk two years before.
Finally,
we discuss the influence on George Griss' negationless
mathematics.
\end{abstract}

\noindent
keywords:
Luitzen Egbertus Jan Brouwer,
George Griss,
intuitionism,
Jules Molk,
principle of the excluded middle

\chapter{Rationale for this translation}

In his seminal paper 
`The unreliability of the logical principles'
\citeyearpar{Brouwer1908C},
Brouwer draws for the first time the revisionistic
consequences of the general view on logic that he had
presented in his
dissertation
\citep[pp.~125–132]{Brouwer1907},
by rejecting the principle
of excluded middle.
The paper appeared in Dutch;
Brouwer's first
published remarks in more widely read languages on the
unreliability
of the principle of excluded middle occur in
\citet[p.~92n2, p.~96n1]{Brouwer1913A}
in English,
and in
\citet[p.~80]{Brouwer1914}
in German.%
\footnoteA{The Brouwer
archive contains an as yet unpublished German translation of the paper at hand,
most likely
prepared by Karl Menger in 1925 or 1926,
who in those
two years was Brouwer's assistant.
In `My memories of L.E.J.~Brouwer',
he writes in footnote 11:
`Brouwer's very moderate assignments to me
were essentially confined to translations of some older
writings of his
on intuitionism from Dutch into German.'
\citep[p.~252]{Menger1979}.}

An English translation,
by Heyting and Gibson,
appeared in 1975,
in volume 1 of Brouwer's
\emph{Collected Works} \citep[pp.~107–111]{Brouwer1975}.
The project of attempting a novel translation
seemed to us a worthwhile one
against the following background.
Brouwer's original text must have struck already a Dutch reader in
1908
as a difficult and unusual one,
whose author nevertheless retains a full mastery of his
sentences.
In order to preserve
for a reader of English at least part of
what the original text thus
conveys to a reader of Dutch,
we believe that one must,
to put it in Schleiermacher's memorable terms,
move the reader towards the author,
instead of moving the author towards the reader.
We have therefore aimed to translate as literally as
possible;
to translate Dutch words by English cognates,
and to preserve relations between Dutch cognates among their
English translations,
wherever appropriate;
to preserve the Germanic structure of the original
Dutch to the extent that English,
likewise a Germanic language,
allows for it;
and to preserve Brouwer's idiomatic idiosyncracies.

\chapter{Brouwer's submission of his manuscript}

From its novel treatment of the principle of
excluded middle,
it is clear that Brouwer drafted his manuscript after
the thesis \citep{Brouwer1907},
which was defended on February 19,
1907;
and towards the
end of that year,
he submitted it to the
\emph{\dutch{Tijdschrift voor Wijsbegeerte}}.%
\footnoteA{From 1933,
that journal appeared as
\emph{\dutch{Algemeen Nederlands Tijdschrift voor Wijsbegeerte en Psychologie}},
and from
1970 onwards as
\emph{\dutch{Algemeen Nederlands Tijdschrift voor Wijsbegeerte}}.
It must not be confused with the
\emph{\dutch{Tijdschrift voor Filosofie}}
that has been published in Leuven since
1939.}
In
a letter to one of its editors,
the philosophically inclined man of
letters Johannes Diderik Bierens de Haan,
of December 7,
Brouwer had promised to explain matters further in
subsequent papers that would be longer and better
understandable to non-mathematicians.
That letter has not survived,
but
this element of it is taken up in another letter that did.
On January 3,
1908,\footnoteA{\citet[Online Supplement, pp.~225–226]{Dalen2011}.}
another editor,
the physicist P.~Kohnstamm,
informed Brouwer that the paper had been accepted that day,
in spite of most editors confessing to have understood very
little of
it.
Van Dalen
\citeyearpar[p.~108]{Dalen1999}
suspects that the editorial board had also
hesitated to publish Brouwer's paper because he was not a
professional
philosopher.
Moreover,
the Dutch professional philosophers had not
appreciated the attempt by the student Brouwer,
two years earlier,
to found a philosophy journal together with Mannoury.
To complicate
matters further,
the \emph{\dutch{Tijdschrift voor Wijsbegeerte}} was
the very journal that was founded in reaction to their
initiative.
On the other hand, it seems that any storm there might have been
had blown over soon,
as the title page of the first volume of the \emph{\dutch{Tijdschrift voor Wijsbegeerte}}
lists Brouwer and Mannoury among the people
`who have promised to contribute'.

In his letter,
Kohnstamm added that he had succeeded in making the
case for acceptance mainly because of Brouwer's promise to Bierens
de Haan. It is not clear whether Brouwer ever undertook to write the
projected sequels. Also, Kohnstamm gave Brouwer the option of adding
elucidations to the accepted manuscript; for lack of relevant archive
material, we cannot tell whether the published version differs from the
manuscript originally submitted.

Kohnstamm had just published a criticism of psychologism in logic
in the
\emph{\dutch{Tijdschrift}} \citep{Kohnstamm1907},
in the form of a negative review of
Gerard Heymans'
\emph{\german{Die Gesetze und Elemente des wissenschaftlichen Denkens}}
\citep{Heymans1890-1894}.%
\footnoteA{Among the references in
\citet{Kohnstamm1907}
we note on pp.~387–388
Husserl's \emph{\german{Logische Untersuchungen}}
and
Meinong's \emph{\german{Gegenstandstheorie}},
and,
indeed,
on p.~403,
Brouwer's dissertation~–
`a stern Dutch book'
(\dutch{\emph{een ernstig Hollandsch boek}}).
However,
that reference just concerns Brouwer's adoption
\citep*[p.~104]{Brouwer1907}
of Poincaré's remark that the world started rotating only
with Copernicus.}
In a letter  of January 18, 1907 to his thesis adviser,
Diederik Johannes Korteweg,
Brouwer had likewise expressed an anti-psychologistic stance:
\begin{quoting}
From your characterization of theoretical logic as part of psychology
I gathered that I had expressed myself rather vaguely,
because it was actually my intention to show
that theoretical logic on no account has a psychological meaning,
even though it is a science. \citep[p.~37, trl.~Van Dalen]{Dalen2011}
\end{quoting}
Brouwer did not discuss the matter in the `Unreliability' paper,%
\footnoteA{Much later, in 1949, Brouwer mentioned in a letter to David van Dantzig
his `belief that psychological pictures of intuitionistic mathematics, however interesting they may be, never can be adequate' \citep[p.~439]{Dalen2011}.
For anti-psychologistic readings of Brouwer,
see \citet{Placek1999} and \citet[ch.~6]{Atten2004a}.}
which is the more regrettable since
the issue of the
\emph{\dutch{Tijdschrift}}
in which it appears contains also
Heymans' reply to Kohnstamm.%
\footnoteA{\citet{Heymans1908}, where Brouwer, Husserl, and Meinong are not mentioned.}

\chapter{Brouwer's conception of logic}

The conception of logic involved in Brouwer's remarks
on the use of the principle of the excluded middle in mathematics
is that formulated in his dissertation.
Logic,
according to Brouwer,
is the study of patterns in linguistic records of
mathematical
acts of construction,
and,
as such,
a form of applied mathematics.
Mathematical constructions out of the
intuition of time
are themselves not of a linguistic nature.
Language cannot play a creative role in mathematics;
there are no mathematical truths that can be arrived at
by linguistic means
(such as logic)
that could not,
at least in principle,
have been arrived at in acts of languageless mathematical
construction \citep[p.~133]{Brouwer1907}.

A correct inference is one where the construction required
by its conclusion can be found from hypothetical actual
constructions
for its premisses.
The hypotheses here are epistemic ones,
in that the
premisses are \emph{known}.
Thus,
they differ from assumptions of the
usual natural deduction kind,
which merely assume that propositions are
true.
For Brouwer's conception of truth,
however,
only these epistemic
assumptions play a role,
since for him to assume that a proposition is
true is to assume that one has a demonstration of it,
that is,
that one
knows that it is true.\footnoteA{For a discussion,
see \citet[sections 1 to 4]{Sundholm.Atten2008}.}

Our use of logical signs
in what follows is meant only as an abbreviatory device.
Although Brouwer in his dissertation had remarked
of the language accompanying logical reasonings that
`As well as any mathematical language this language can without 
much trouble be condensed
into symbols' \citep[p.~159]{Brouwer1975},%
\footnoteA{As an example,
he refers to \citet[p.~35ff]{Whitehead1898}. `\dutch{Zoo goed als alle wiskundige taal is ook
deze taal zonder moeite te condenseeren tot symbolen}' 
\citep[p.~159n]{Brouwer1907}.}
in his own writings he persisted in preferring 
sometimes prolix non-symbolic language.
We will write
$A \rightarrow B$
for
`A 
(hypothetical)
actual
construction for $A$
can be continued into
a construction for $B$'.%
\footnoteA{Heyting's later works on logic do employ symbolism,
but here we leave it open
whether they are committed to this meaning.
See for example \citet{Sundholm1983}.}

The logical principles referred to in
the title of Brouwer's paper are
those of
the Aristotelian tradition:
the principles of the syllogism
(in the
paper defined by modus Barbara),
of contradiction,
and of the excluded
third.\footnoteA{Brouwer's choice of terminology here is
different from
that of his likely logic teacher at Amsterdam,
C.F.~Bellaar-Spruyt.
The latter's posthumously published book on formal logic
\citep{Bellaar-Spruyt1903}
lists the principium identitatis,
principium
contradictionis,
principium exclusi medii
(also `tertium vel medium non
datur',
\citealt[p.~18]{Bellaar-Spruyt1903}),
and the dictum de omni et nullo,
\citealt[p.~14]{Bellaar-Spruyt1903}.
Brouwer's principle
of the syllogism seems to comprise the principium
identitatis and the
dictum de omni et nullo.
The principles of identity and of syllogism are
also discussed by Poincaré in `\french{Sur la nature du
raisonnement
mathématique}' \citep{Poincare1894},
which also appears,
in abridged form,
in the first
chapter of
\emph{\french{Science et Hypothèse}} \citep{Poincare1902};
Brouwer knew in
any case the latter of the two.}
Of course,
from lectures by Gerrit Mannoury
Brouwer knew about further developments,
in particular those by Frege
and Peano;\footnoteA{These lectures were published in
shortened and
revised form as \emph{\german{Methodologisches und
Philosophisches zur
Elementar-Mathematik}} \citep{Mannoury1909}.}
but for his principled
criticism it suffices to consider the Aristotelian case.

\chapter{Unreliability in the natural sciences and in wisdom}

Brouwer introduces the main question of his paper,
that of the reliability
of logic in pure mathematics, by arguing that in two other domains logic
is not reliable: the natural sciences and wisdom.

The problem with the use of logic in the natural sciences, as Brouwer
describes it, is the familiar problem of induction. There is no guarantee
that a mathematical model that explains a given set of observations will
correctly predict further observations. But logic leads from statements
in the mathematical model to other statements in that model. Hence, it
may well lead from premisses that agree with observations to conclusions
that do not, and is, in that sense, unreliable.

In wisdom, logic is not reliable for a different type of reason. Logic
presupposes the presence of mathematical constructions, but in wisdom such
constructions are absent. Mathematics embraces time awareness, whereas
wisdom discards it.\footnoteA{Compare this remark in
`Will, knowledge and speech':\csname @beginparpenalty\endcsname10000
\begin{quoting}
Mathematical attention is not a necessity
but a phenomenon of life subject to the free will,
everyone can find this out
for himself by internal experience:
every human being can at will either dream-away time-awareness
and the separation between the Self and the World-of-perception
or
by his own powers
bring about this separation
and
call into being in the world-of-perception the condensation of separate things.
\citep[pp.~418–419]{Stigt1990}
\end{quoting}
(\dutch{`Dat wiskundige beschouwing geen
noodzaak, doch een aan den vrijen wil onderworpen levensverschijnsel is,
daarvan kan ieder bij zichzelf de inwendige ervaring opdoen:
ieder mensch kan naar willekeur hetzij zich zonder tijdsgewaarwording
en zonder scheiding tusschen Ik en Aanschouwingswereld verdroomen,
hetzij de genoemde scheiding door eigen kracht voltrekken en in de
aanschouwingswereld de condensatie van aparte dingen in het leven
roepen'}, \citealt[p.~2]{Brouwer1933A2}.)}
Since time awareness introduces the subject-object
distinction, it keeps consciousness out of what Brouwer later called
its `deepest home' \citep[p.~1235]{Brouwer1949C}. An attempt to apply logic
to wisdom would require one to impose a mathematical structure on it,
thereby distorting its content. Logic is unreliable in this domain,
for logical conclusions from distorted content cannot be expected to
reflect that content accurately.%
\footnoteA{On Brouwer's interest
in religion,
mysticism,
and their relations to science,
see \citet[sections 1.3 and 1.6]{Dalen1999};
\citet{Stigt1996};
and \citet{Koetsier2005}.
A comparison with Gödel on this point is presented in
\citet{Atten.Tragesser2003}.}

The question then arises whether in pure mathematics, where, in contrast
to natural science, abstraction has been made from all observational
content, and, in contrast to wisdom, logic is applied to something that
does have mathematical structure, the use of logic is reliable. The main
point of this paper is that it is not.

\chapter{Unreliability in mathematics}

Brouwer had already made a case for the possibility of unreliable logical principles in
his dissertation:
\begin{quoting}
And if one succeeds in the construction of
\emph{linguistic}
buildings,
sequences of sentences proceeding according to the logical
laws,
thereby departing from linguistic images which could
accompany basic mathematical truths in actual mathematical
buildings,
and if it turns out that those linguistic buildings can
never produce the linguistic form of a contradiction,
then all the same they belong to mathematics only in their
quality of a linguistic building,
and have nothing to do with mathematics outside of that
building,
e.g.~with ordinary arithmetic or geometry.

So the idea that by means of such linguistic buildings we
can obtain any knowledge of mathematics apart from that
which can be constructed directly on the basis of intuition,
is mistaken.
And more so is the idea that in \emph{this} way we can lay
the \emph{foundations} of mathematics,
in other words that we can ensure the reliability of the mathematical theorems.
\citep[pp.~132–133, original emphasis]{Brouwer1975}%
\footnoteA{%\begin{quoting}\vspace{-\baselineskip}
\dutch{En wanneer het gelukt
\emph{taal}gebouwen op te trekken,
reeksen van volzinnen,
die volgens de wetten der logica op elkaar volgen,
uitgaande van taalbeelden,
die voor werkelijke wiskundige gebouwen,
wiskundige grondwaarheden zouden kunnen accompagneeren,
en het blijkt dat die taalgebouwen nooit het taalbeeld 
van een contradictie zullen kunnen vertoonen,
dan zijn ze toch alleen wiskunde als taalgebouw
en hebben met wiskunde buiten dat gebouw,
bijv.~met de gewone rekenkunde of meetkunde niets te maken.}

\hspace{\footparindent}\dutch{Dus in geen geval mag men denken,
door middel van die taalgebouwen iets van andere wiskunde,
dan die direct intuitief op te bouwen is,
te kunnen te weten komen.
En nog veel minder mag men meenen,
op 
\emph{die}
manier de
\emph{grondslagen}
der wiskunde te kunnen leggen,
m.a.w.~de betrouwbaarheid der
wiskundige eigenschappen te kunnen verzekeren.}
\citep[pp.~132–133]{Brouwer1907}%\vspace{-\baselineskip}
%\end{quoting}
}
\end{quoting}
And in entry \textsc{xx} in the list of propositions
submitted to the public defence together with it,
according to Dutch custom that is still today 
observed in some universities,
he had said:
\begin{quoting}
To secure the reliability of mathematical reasonings
one cannot succeed solely
by starting from some sharply formulated axioms
and further strictly adhering to the laws of theoretical logic.
\citep[p.~101]{Brouwer1975}%
\footnoteA{\dutch{Het kan niet gelukken,
de betrouwbaarheid der wiskundige redeneeringen te verzekeren,
enkel door uit te gaan van eenige scherp gestelde axioma's
en verder streng vast te houden aan de wetten
der theoretische logica}. \citep[\dutch{Stellingen}]{Brouwer1907}}
\end{quoting}
The
reliability of logical reasoning depends on the mathematical
context
in which it is applied:
it is the context that determines whether
the logical reasonings can be traded in for corresponding
mathematical
constructions.
In the dissertation Brouwer rejected the attempt to come to
know,
by the use of logic,
something mathematical that is nonconstructive;
he there considered
reliable within constructive
mathematics not only the principles of the syllogism
\citep[p.~131]{Brouwer1907}
and of
contradiction, but also the principle of excluded middle.
The reason
is that at the time he read $A \vee \neg A$ as $\neg A
\rightarrow
\neg A$:%
\footnoteA{As Van Dalen \citeyearpar[pp.~106–107]{Dalen1999} has pointed out,
Brouwer most likely arrived at this reading under the influence
of the logic lectures by Bellaar-Spruyt.}
\begin{quoting}
While in the syllogism a mathematical element could be discerned,
the proposition:

A function is differentiable or is not differentiable

\noindent
says
\emph{nothing};
it expresses the same as the following:

If a function is not differentiable,
then it is not differentiable.

But the logician,
looking at the
\emph{words}
of the former sentence,
and discovering a regularity in the combination of words
in this and in similar sentences,
here again projects a mathematical system,
and he calls such a sentence an
\emph{application of the tertium non datur}.
\citep[p.~75, original emphasis]{Brouwer1975}%
\footnoteA{\dutch{Was in het syllogisme nog een wiskundig element te
onderkennen,
de stelling:}

\hspace{3ex}\dutch{Een functie is òf differentieerbaar òf niet
differentieerbaar}

\dutch{zegt 
\emph{niets};
drukt hetzelfde uit,
als het volgende:}

\hspace{3ex}\dutch{Als een functie niet differentieerbaar is,
is ze niet differentieerbaar.}

\dutch{Maar de 
\emph{woorden}
van eerstgenoemde volzin bekijkend, en een
regelmatig gedrag in de opvolging der woorden van deze en
van dergelijke volzinnen ontdekkend, projecteert
de logicus ook hier een wiskundig systeem, en noemt zulk een
volzin een 
\emph{toepassing van het principe van tertium nondatur}.}
\citep[p.~131]{Brouwer1907}}
\end{quoting}

\chapter{The principle of excluded middle is unreliable}

In `Unreliability',
Brouwer will advance upon the dissertation in two ways:
he corrects his reading of the principle of excluded middle,
and he
shows that this corrected understanding entails the unreliability of a
traditional principle within constructive mathematics itself.
% see also Brouwer 1949C (p.~1243).

This is the
(silently)
corrected understanding of the principle of excluded middle:
\begin{quoting}
Now the principium \emph{tertii exclusi}: this demands that every
supposition\footnoteA{Suppositions should here not be taken
in the sense of abstract propositions in a Platonic realm of
abstract entitites, as in Bolzano or in Frege. What seems to be
meant is rather: Every mathematical assumption that we can make
is either correct or incorrect.} is either correct or incorrect,
mathematically: that of every supposed fitting in a certain way
of systems in one another, either the termination or the blockage
by impossibility, can be constructed.
\citep[p.~156, trl.~ours]{Brouwer1908C}%
\footnoteA{\dutch{Nu het principium \emph{tertii exclusi}:
dit eischt,
dat iedere onderstelling òf juist òf onjuist is,
wiskundig:
dat van iedere onderstelde inpassing van systemen op
bepaalde wijze in elkaar hetzij de beëindiging,
hetzij de stuiting op onmogelijkheid kan worden
geconstrueerd.}}
\end{quoting}
The change of
mind is acknowledged in
`\dutch{Addenda and corrigenda to “On the
Foundations of Mathematics}”~'
\citep[p.~1]{Brouwer1917A2}.

Note that $\neg A$ does not merely mean
that no proof of $A$ exists,
but
that from an assumed actual demonstration of $A$ one can
`construct the blockage by impossibility'
(see also \citealt[p.~127]{Brouwer1907}).
In this sense, intuitionistic negation is unlike the
classical notion a positive notion,
as it involves the existence of a blockage.%
\footnoteA{\citet[pp.~498–500]{Becker1927},
with reference to the section `Evidence and truth'
(\emph{\german{Evidenz und Wahrheit}})
in the sixth of Husserl's
\emph{\german{Logische Untersuchungen}}
\citep{Husserl1984b}.
This is clearly the passage by Becker that Heyting has in mind in
Königsberg
\citep[p.~113]{Heyting1931}:
\begin{quoting}
\german{Eine logische Funktion ist ein Verfahren,
um aus einer gegebenen Aussage eine andere Aussage
zu bilden.
Die Negation ist eine solche Funktion;
ihre Bedeutung hat Becker,
im Anschluß an Husserl,
sehr deutlich beschrieben.
Sie ist nach ihm etwas durchaus Positives,
nämlich die Intention
auf einen mit der ursprünglichen Intention
verbundenen Widerstreit.}
\end{quoting}
(`A logical function is a method for turning a given statement into
another statement.
Negation is such a function;
Becker,
following Husserl,
has described 
its meaning very clearly.
It is according to him something wholly positive,
namely the intention directed to a conflict
bound up with the original intention.' Trl.~ours.)
Heyting does not give a reference here,
but had already mentioned Becker's
\emph{\german{Mathematische Existenz}}
on p.~107.}

Thus
understood,
the principle of excluded middle is not reliable,
for we
do not have a general decision method as required by the
constructive
reading.
Brouwer's claim is not that we can never have such a method:
`in infinite systems the principium tertii exclusi is
\emph{as yet} not
reliable'
(our emphasis).
Brouwer states the first so-called `Brouwerian
counterexamples' or `weak counterexamples' to the principle
of excluded
middle,
which illustrate its unreliability.
These are propositions of
which we are in a position to assert the weak negation,
but not the truth
or the strong negation.
Of course any open problem is,
as such,
a weak
counterexample to the principle of excluded middle;
the importance of
weak counterexamples comes from the fact that they can be
used to show
that certain highly general principles have not yet been
established,
such as `Every set is finite or infinite' or `The continuum
is totally
ordered'.
Brouwer published weak counterexamples to the principle of
excluded middle also in international journals,
but only much later
\citep{Brouwer1921A, Brouwer1924N, Brouwer1925E,
Brouwer1929A}.
By then
he had found a uniform technique for constructing weak
counterexamples
that depended only on the fact \emph{that} open problems of
a certain
simple type still exist,
not on the exact content of these problems.

\chapter{Are there absolutely undecidable propositions?}

Brouwer adds
the following remark to his explanation of the principle of excluded
middle:
\begin{quoting}
The question of the validity of the principium tertii exclusi is
thus equivalent to the question concerning the \emph{possibility
of unsolvable mathematical problems}. For the already proclaimed
conviction that unsolvable mathematical problems do not exist,
no indication of a demonstration is present.\footnoteA{Compare
entry \textsc{xxi} in the list of theses in the dissertation:
`\dutch{Ongegrond is de overtuiging van Hilbert (Gött.~Nachr.~1900, pag.~261):
“\german{dass ein jedes bestimmte mathematische Problem einer strengen Erledigung notwendig fähig sein
müsse, sei es, dass es gelingt, die Beantwortung der
gestellten Frage zu geben, sei es dass die Unmöglichkeit der Lösung und damit die Notwendigkeit des
Misslingens aller Versuche dargetan wird}”}.' \citep[\dutch{Stellingen}]{Brouwer1907}}
\citep[p.~5]{Brouwer1908C}
\end{quoting}
The claim that every mathematical problem is solvable is of
course constructively stronger than the claim that there are no
unsolvable problems.\footnoteA{See also \citet[6.5]{Wittgenstein1922},
\citet{Schlick1935},
\citet{McCarty2005},
and
\citet{Martin-Lof1995} (in particular the postscript in the reprint in \citealt{Schaar2012}).} The former is equivalent to the principle
that for any $A$, $A \vee \neg A$, the latter to the principle that
for any $A$, $\neg\neg(A \vee \neg A)$; and Brouwer had demonstrated
the validity of the latter in the same paper. Indeed, in the Brouwer
archive there is a note from about the same period 1907–1908 in which
the point is made explicitly:
\begin{quoting}
Can one ever demonstrate of a question, that it can never
be decided? No, because one would have to do so by reductio
ad absurdum. So one would have to say: assume that the
question has been decided in sense $a$, and from that deduce
a contradiction. But then it would have been demonstrated that
not-$a$ is true, and the question remains decided. \citep[p.~174n.~a, trl.~ours]{Dalen2001}\footnoteA{\dutch{Zal men nu ooit van een vraag kunnen bewijzen,
dat ze nooit uitgemaakt kan worden?
Neen,
want dat zou moeten uit het ongerijmde.
Men zou dus moeten zeggen:
Gesteld dat het was uitgemaakt in zin $a$
en daaruit afleiden,
tot een contradictie kwam.
Dan zou echter bewezen zijn,
dat niet $a$ waar was,
en de vraag bleef uitgemaakt.}}
\end{quoting}

% file: uitgesloten_midden
Brouwer never published this note. Wavre in 1926 gave the argument for
a particular case, while clearly seeing the general point:
\begin{quoting}
It suffices to give an example of a number of which one does
not know whether it is algebraic or transcendent in order
to give at the same time an example of a number that, until
further information comes in, \emph{could} be neither the one
nor the other. But, on the other hand, it would be in vain, it
seems to me, to want to define a number that \emph{is} neither
algebraic nor transcendent, as the only way to show that it is
not algebraic consists in showing that it is absurd that it would
be, and then the number would be transcendent. \citep[p.~66, trl.~ours, original emphasis]{Wavre1926}\footnoteA{\french{Il suffit
donc de fournir l'exemple d'un nombre dont on ne sache s'il est
algébrique ou transcendant pour fournir en même temps l'exemple
d'un nombre qui, jusqu'à plus ample information, \emph{pourrait}
n'être ni l'un ni l'autre. Mais, d'autre part, il serait vain,
me semble-t-il, de vouloir définir un nombre qui ne \emph{soit}
ni algébrique ni transcendant, car la seule manière de prouver
qu'il n'est pas algébrique consistant à prouver qu'il serait
absurde qu'il le fût, ce nombre serait transcendant.} [original
emphasis]}
\end{quoting}
The general, schematic point was explicitly noted by Heyting in 1934:
\begin{quoting}
Further, the formula $\vdash\neg\neg(a\vee\neg a)$ should
be highlighted. It has the same meaning\footnoteA{Heyting
writes `\german{gleichbedeutend}'. Note that the two formulas are
equi-assertible, but have different assertion conditions.}
as $\vdash\neg(\neg a \wedge \neg\neg a)$ and expresses
Brouwer's theorem on the absurdity of the absurdity of the
excluded third, and amounts to saying that a demonstrably unsolvable
problem cannot exist. \citep[p.~16]{Heyting1934}%
\footnoteA{\german{Es sei noch die
Formel $\vdash\neg\neg(a\vee\neg a)$ hervorgehoben, die mit
$\vdash\neg(\neg a \wedge \neg\neg a)$ gleichbedeutend ist und den
Brouwerschen Satz von der Absurdität der Absurdität des Satzes
vom ausgeschlossenen Dritten zum Ausdruck bringt. Sie besagt,
daß es ein nachweisbar unlösbares Problem nicht geben kann.}}
\end{quoting}

\chapter{The principle of excluded middle is consistent}

Brouwer
also observes that, although the principle of excluded middle is not
schematically valid, none of its instances is false, since $\neg(A \vee
\neg A)$ implies the contradiction $\neg A \wedge \neg\neg A$. This
demonstrates the correctness of the principle that, for any $A$,
$\neg\neg(A \vee \neg A)$. Brouwer concludes that it is always consistent
to use (this form of) the principle of excluded middle but that it does
not always lead to truths. Later, Brouwer gave a refutation of the schema
$\forall x(P(x) \vee \neg P(x))$ using specifically intuitionistic
principles regarding choice sequences and continuity \citep{Brouwer1928A2}.\footnoteA{In the Bishop tradition,
some versions of the principle
of excluded middle that Brouwer devised counterexamples to have been
given a systematic place:
\textsc{lpo},
\textsc{wlpo},
\textsc{llpo}\@.
See \citet[ch.~1, section 1]{Bridges.Richman1987}.}
Hence Brouwer's proposal to divide the theorems that are usually
considered as having been demonstrated into the correct and the
non-contradictory ones \citep[7n.~2]{Brouwer1908C}, that is, those whose
reduction to absurdity has been refuted. That is not a suggestion
that there are three truth values, true, non-contradictory,
false;\footnoteA{Church makes this point clearly in 1928 in `On the
law of excluded middle' \citep{Church1928}, criticising \citet{Barzin.Errera1927}.} for a non-contradictory proposition might be proved
one day and thereby become true. This observation on the consistency of
the principle of excluded middle would, in the 1920s, be at the basis
of Brouwer's optimism,
expressed in print but nevertheless often neglected,
concerning the success of the Hilbert Program,
a success that Brouwer would consider of no general mathematical value
(\citealt[p.~3]{Brouwer1924N}; \citealt[p.~377]{Brouwer1928A2}).

\chapter{Partial validity of the principle of excluded middle}

The principle of excluded middle is valid, Brouwer points out, in
finite domains, for questions whether a given construction of finite
character\footnoteA{Not, obviously, in the sense of the classical equivalent of
Zermelo's axiom of choice known as the Teichmüller-Tukey lemma.} is
possible. Only finitely many attempts at that construction can be made,
and each will succeed or fail in finitely many steps (see also \citealt[p.~114]{Brouwer1955}). Brouwer came to explain the genesis of the belief in the
validity of the principle of excluded middle as follows:
\begin{quoting}
I am convinced that the axiom of solvability and the principle
of excluded third are both false,\footnoteA{Brouwer came to
identify those, but after they had been presented separately.}
and that historically the belief in these dogmas has been caused
thusly. First, one has abstracted classical logic from the
mathematics of subsets in a certain finite set, then ascribed
to this logic an a priori existence independent of mathematics,
and finally, on the basis of this alleged apriority, applied
it rightlessly to the mathematics of infinite sets.
\citep[n.~4]{Brouwer1922A}\footnoteA{\german{Meiner Ueberzeugung nach
sind das Lösbarkeitsaxiom und der Satz vom ausgeschlossenen
Dritten beide falsch und ist der Glaube an diese Dogmen historisch
dadurch verursacht worden, dass man zunächst aus der Mathematik
der Teilmengen einer bestimmten endlichen Menge die klassische
Logik abstrahiert, sodann dieser Logik eine von der Mathematik
unabhängige Existenz a priori zugeschrieben und sie schliesslich
auf Grund dieser vermeintlichen Apriorität unberechtigterweise
auf die Mathematik der unendlichen Mengen angewandt hat.}}
\end{quoting}
% to an unwarranted projection from such finite cases (in particular,
% those arising from the application of finite mathematics to everyday
% phenomena) to the infinite.
See on this point also Brouwer (\citeyear[p.~2]{Brouwer1924N}; \citeyear[pp.~423–424]{Brouwer1929A};
\citeyear[p.~492]{Brouwer1949C}; and \citeyear[pp.~510–511]{Brouwer1952B}).

\chapter{Brouwer's concern with meaning and truth}

% check also if McCarty (eg Cerisy) makes remarks similar to Koss'

We wish to discuss one further aspect of the text itself.
In it,
Brouwer uses neither
the term
`\dutch{(wiskundige) waarheid}'
((mathematical) truth),
nor
`\dutch{betekenis}'
(meaning);
but
a careful consideration of his Dutch
and some of his other writings
will reveal that these really are the notions under consideration.

We begin with truth.
In place of
`\dutch{waar}' (true) and
`\dutch{onwaar}' (false),
Brouwer uses
`\dutch{juist}' and
`\dutch{onjuist}'.
The term
`\dutch{juist}',
however,
is most commonly translated by `right' and/or `correct',
which raises the question as to how this relates to (propositional) truth.
The largest and most authorative dictionary
of the Dutch language, 
the \emph{\dutch{Woordenboek der Nederlandsche Taal}},
lists among the meanings of `\dutch{juist}':
`\dutch{Met de waarheid —,
met
het wezen van iets in overeenstemming;
de waarheid weergevende;
aan
de waarheid beantwoordende}'
(in agreement with the truth,
with the essence of something;
representing the truth;
corresponding to the truth),
and gives among its historical examples 
this sentence from 1897:
`\dutch{Met “waarheid” kan men bedoelen de meest
juiste voorstelling van de dingen}'
(By `truth' one can mean the most correct representation of things)~–
by a happenstance,
written by an author who would become one of
Brouwer's best friends,
Frederik van Eeden.
Here,
the sense of
`\dutch{juist}'
is given by the translation `true',
chosen by Heyting in the \emph{Collected Works}
\citep[p.~110]{Brouwer1975}.
We prefer the alternative `correct',
however,
in order to translate
the different words
`\dutch{juist}'
and
`\dutch{waar}'
differently.

Of course correctness can be relative to something else than truth,
for example a convention,
a value,
or an ideal;
but that here truth is meant
is clear from the fact that
in his dissertation Brouwer indeed is willing to speak of
`\dutch{wiskundige grondwaarheden}'
(basic mathematical truths),%
\footnoteA{\citet[p.~132]{Brouwer1907}.}
in his reply to Mannoury's review of his dissertation of
`\dutch{wiskundige waarheden}'
(mathematical truths),%
\footnoteA{\citet[p.~328]{Brouwer1908B}.}
and
in his draft letter to De Vries,
dated February 15, 1907,
of 
`\dutch{de waarheid van de wiskundige stellingen}'
(the truth of the mathematical theorems).%
\footnoteA{\citet[Online Supplement, p.~201]{Dalen2011}.}
An explicit identification is found in
`\dutch{Willen, weten, spreken}'
(Will, knowledge and speech) of 1932:
`\dutch{juiste
(d.w.z.~daadwerkelijk wiskundige beschouwingen doeltreffend indicerende)
affirmaties}'
\citep[p.~54]{Brouwer1933A2},
translated by Van Stigt%
\footnoteA{The \emph{Collected Works} do not give an English
translation of the first two sections of `Will, knowledge and speech'
because of their proximity to the corresponding sections of the first Vienna lecture \citep{Brouwer1929A},
which they include in German.}
as
`correct affirmations
(i.e.~effectively indicating actual mathematical viewing)'
\citep[p.~424]{Stigt1990},%
\footnoteA{At the corresponding place in the first Vienna lecture,
Brouwer had written:
`\german{zutreffenden
(d.h.~tatsächliche mathematische Betrachtungen andeutenden)
Aussagen}' \citep[p.~158]{Brouwer1929A}.
We here note that
the range of meaning of
the German `\german{zutreffend}' is included in and
much narrower than that of the Dutch `\dutch{juist}'.}
where the viewing takes place in languageless intuition.%
\footnoteA{The visual metaphor is not often used in Brouwer's writings.
In the dissertation one finds:
\begin{quoting}
Now we have 
seen that classical logic studies the linguistic counterpart
of logical reasonings, 
i.e.~of reasonings on
\emph{relations of whole and part}
for arbitrary mathematically constructed systems;
from the fact that we
\emph{see}
these mathematical systems we may 
conclude that
\emph{here}
the sentences succeeding one another
according to classical 
logic, will never show contradictions,
because they correspond to acts of mathematical construction.
\citep[p.~88, original emphasis]{Brouwer1975}
\end{quoting}
(`\dutch{Nu hebben we gezien,
dat de klassieke logica bestudeert de taalbegeleiding
der logische redeneeringen,
d.w.z.~der redeneeringen
\emph{in relaties van geheel en deel}
voor willekeurige wiskundig opgebouwde systemen;
en we weten uit het feit,
dat we die wiskundige systemen
\emph{zien},
dat
\emph{daar}
de volgens de klassieke logica
elkaar opvolgende volzinnen,
die immers wiskundige bouwhandelingen begeleiden,
nooit contradicties zullen vertoonen}',
\citealt[pp.~159–160, original emphasis]{Brouwer1907}.)}

We now turn to meaning.
As we have seen,
Brouwer begins his discussion of the principle of the excluded middle
as follows:
\begin{quoting}
Now the principium \emph{tertii exclusi}:
this demands
that every
supposition is either correct or incorrect…
\end{quoting}
In the original:
\begin{quoting}
Nu het principium \emph{tertii exclusi}:
dit eischt,
dat iedere onderstelling òf juist òf onjuist is…
\end{quoting}
Brouwer writes `demands' (\dutch{\emph{eischt}}), 
not,
as one might have expected,
`means' (\dutch{\emph{betekent}}) or `asserts' (\dutch{\emph{beweert}}).
However,
the
\emph{Woordenboek der Nederlandsche Taal}
lists
among the meanings of `\dutch{eischen}':
`\dutch{Tot voorwaarde voor zijn bestaan,
welvaren of welslagen hebben}':
`to have as a condition for its
existence,
prospering,
or
success'.
Hence,
`\dutch{dit principe eischt dat}…' can naturally
be understood as
`for this principle to hold,
what is
required is that…'.
But this condition clearly amounts to a meaning specification of the principle.

That Brouwer intends this sense of
`\dutch{eischt}'
is brought out by a comparison with his discussion of his first principle,
to wit syllogism,
where at the corrresponding place he uses
`\dutch{leest in}'
(reads~… as~…).
The latter
unambiguously expresses a concern with meaning.
A coherent
interpretation of Brouwer's remarks
should
accord the same sense to these two verbs;%
\footnoteA{Brouwer's remark on the second principle does not
contain a corresponding verb at all.}
the disambiguation of
`\dutch{eischt dat}' should pick out the sense in which it has
the same meaning
as `\dutch{leest in}'.%
\footnoteA{\citet[p.~75]{Koss2013}
reaches the conclusion
that in this paper Brouwer's concern was neither with meaning nor
with truth.
As we show above,
the linguistic facts do not bear him out.}

% something about meaning in the dissertation,
% and at the thesis defence (`zinloze figuurtjes of gestamel')?

\chapter{Precursors}

Brouwer was not the first who voiced criticism or
hesitations about either the usefulness or the validity of
the principle
of excluded middle in a purely mathematical context.
From the 1870s,
Kronecker objected to the unlimited use of the principle of
excluded
middle and of definition by undecided separation of cases.
For example,
in
his treatise on algebraic numbers of 1882,
he wrote on the factorization
of polynomial functions:
\begin{quoting}
The definition of irreducibility drawn up in section 1 lacks
a secure grounding as long as no method has been indicated by
which it can be decided whether a definite given function is
irreducible according to that definition or not.
\citep[pp.~10–11]{Kronecker1882}\footnoteA{\german{Die im Artikel 1 aufgestellte Definition
der Irreduktibilität entbehrt solange einer sicheren Grundlage,
als nicht eine Methode angegeben ist, mittels deren bei einer
bestimmten vorgelegten Funktion entschieden werden kann, ob
dieselbe der aufgestellten Definition gemäß irreduktibel ist
oder nicht.}}
\end{quoting}
adding in a footnote,
\begin{quoting}
The analogous need, which as a matter of fact has often
remained neglected, arises in many other cases, in definitions
as in demonstrations, and on another occasion I will come
back to this generally and thoroughly.
\citep[p.~11n]{Kronecker1882}\footnoteA{\german{Das analoge Bedürfnis, welches freilich
häufig unbeachtet geblieben ist, zeigt sich in vielen anderen
Fällen, bei Definitionen wie bei Beweisführungen und ich werde
bei einer anderen Gelegenheit in allgemeiner und eingehender
Weise darauf zurückkommen.}}
\end{quoting}
His student Jules Molk\footnoteA{December 8, 1857, Strasbourg~–~May 7,
1914, Nancy.} gave voice to the doubts of his \german{Doktorvater} in the printed
version of his Berlin dissertation from 1885:
\begin{quoting}
The definitions should be algebraic and not only logical. It does
not suffice to say: `A thing exists or it does not exist'. One
has to show what being and not being mean, in the particular
domain in which we are moving. Only thus do we make a step
forward.
\citep[p.~8, trl.~ours]{Molk1885}%
\footnoteA{\french{Les définitions devront être algébriques et non pas
logiques seulement. Il ne suffit pas de dire: `Une chose est
ou elle n'est pas'. Il faut montrer ce que veut dire être et
ne pas être, dans le domaine particulier dans lequel nous nous
mouvons. Alors seulement nous faisons un pas en avant.}}
\end{quoting}
Molk became professor in Nancy, and was the editor-in-chief and driving
force behind the French version of Felix Klein's \emph{Enzyklopädie
der mathematischen Wissenschaften und ihren Grenzgebiete}. He translated
and augmented, especially concerning foundational matters, Pringsheim's
beautiful surveys of topics in elementary analysis. In Book I, volume I.3,
section 10, `\french{Point de vue de L.~Kronecker}', Molk considerably
elaborated upon the above brief remark from his dissertation:
\begin{quoting}
Analysis should, on the other hand, refrain from general
considerations of a logical kind alien to its object. In Analysis
definitions may introduce nothing but auxiliary notions that
facilitate the study of the various natural groups that one forms
to study the properties of numbers. These auxiliary notions
must have an arithmetical character and not a merely logical
one, whence they can only be about \emph{groups of which each
element can be effectively obtained by means of a finite number
of operations}, and not about groups simply determined by a
non-contradictory logical convention.

Similarly, the logical evidentness\footnoteA{The French here has
`\french{évidence}', in the sense of Cartesian clarity. This meaning,
`evidence of',
is the first given by the \textsc{oed}.
The more familiar,
and in Anglo-American philosophy all-pervasive `evidence for'
exists only in English.
For further discussion,
see \citet{Sundholm2014}.}
of a
reasoning does not suffice to legitimize the use of that reasoning
in Analysis. In order to give a mathematical demonstration of a
proposition, it does not suffice, for example, to establish that
the contrary proposition implies a contradiction. One has to give
a procedure that, operating on the elements under consideration,
by means of a finite number of arithmetical operations in
the old sense of the word, permits one to obtain the result
formulated by the proposition to be demonstrated. \emph{This
procedure constitutes the essence of the demonstration; it is
not an addition to it.}

…

The principle of economy in science~– economy of time,
economy of efforts~– in Analysis leads us to the absolute
and relative rational numbers: that introduction is legitimate,
because its only effect is to shorten the deductions without
changing their character. To every proposition about rational
numbers, for example expressed by an equation, corresponds a
congruence taken according to an easily determined module or
system of modules.

…

The character of demonstrations is, on the contrary,
completely changed by the introduction of \emph{arbitrary}
irrational numbers. One cannot, moreover, give any definition
of these numbers except a logical one, determining them,
but \emph{not mathematically defining them}. It is that
logical (but not mathematical) definition that confers on
(infinite) sets of rational numbers that define so-called
arbitrary irrational numbers, the character of an organic
sequence. Therefore, those numbers, according to L.~Kronecker,
cannot rightfully occur in the definitive demonstration of a
proposition of Analysis.
\citep[pp.~159–61, trl.~ours, original emphasis]{Molk1904}\footnoteA{%\vspace{-\baselineskip}
%\begin{quoting}\footnotesize
\french{L'Analyse doit, d'autre part, se garder de
considérations générales d'ordre logique étrangères
à son objet. Les définitions ne doivent introduire
en Analyse que des notions auxiliaires facilitant
l'étude des divers groups naturels que l'on forme
pour étudier les propriétés des nombres. Ces
notions auxiliaires doivent avoir un caractère
arithmétique et non logique seulement, en sorte
qu'elles ne sauraient porter que sur des groups dont
chaque élément puisse être effectivement obtenu
au moyen d'un nombre fini d'opérations, et non sur
des groupes simplement détermines par une convention
logique non-contradictoire.}

\hspace{\footparindent}\french{De même l'évidence
logique d'un raisonnement ne suffit pas pour légitimer
l'emploi de ce raisonnement en Analyse. Pour avoir donné
une démonstration mathématique d'une proposition,
il ne suffit pas, par exemple d'avoir établie que la
proposition contraire implique contradiction. Il faut
donner un procédé permettant d'obtenir, au moyen
d'un nombre fini d'opérations arithmétiques au sens
ancien du mot, effectuées sur les éléments que l'on
envisage, le résultât qu'énonce la proposition à
démontrer. C'est ce procédé qui constitue l'essence
de la démonstration; il ne vient pas s'y ajouter.}

…

\hspace{\footparindent}\french{Le principe de l'économie dans la science,
économie de temps, économie d'efforts, nous amène en
Analyse les nombres rationnels absolus et relatifs:
cette introduction est légitime, puis-quelle n'a pour
effet que d'abréger les déductions sans en changer
le caractère. A toute proposition concernant les
nombres rationnels, exprimée par une égalité par
exemple, correspond une congruence prise suivant un
module ou un système de modules toujours faciles à
déterminer.}

…

\hspace{\footparindent}\french{Le caractère des démonstrations est, au
contraire, complètement changé par l'introduction
de nombres irrationnels quelconques. On ne peut
d'ailleurs donner de ces nombres qu'une définition
logique, les déterminant, mai ne les définissant pas
mathématiquement. C'est cette définition logique
(mais non mathématique) qui confère aux ensembles
(infinis) de nombres rationnels, que l'on dit définir
des nombres irrationnels quelconques, le caractère
d'une suite organique. Ces nombres ne peuvent donc,
suivant L.~Kronecker, figurer a aucun titre dans la
démonstration définitive d'une proposition d'Analyse.}%\vspace{-\baselineskip}
%\end{quoting}
}
%}% weg

% end footnote
\end{quoting}
Any reader of Brouwer's `Unreliability' will
be struck by the coincidence of the views expressed, even down to some
of the finer details: the rejection of indirect existence proofs; the
prohibition of blind, merely symbolic reasoning; the explicit separation
between demonstrated propositions and non-contradictory ones.
In view of this, the question arises whether Brouwer was aware
of Molk's treatment. According to the central Dutch library
catalogue,\footnoteA{\dutch{Nederlandse Centrale Catalogus}.} in
the Netherlands one copy of the original edition of the article by
Pringsheim and Molk was present, namely at the library of the University
of Amsterdam. Unfortunately, that library has not been able to answer
our question exactly when this fascicule of the \french{Encyclopédie} became
available.\footnoteA{Brouwer occasionally made similar inquiries.
Brouwer to his University Library, August 3, 1929:
\begin{quoting}[indentfirst=true]
Undersigned would much appreciate it
if he could learn the exact date
on which Heft 2 of Band 142 of the
\german{Journal für die reine und angewandte
Mathematik} (Crelle’s Journal)	
was received at the University Library.
The date probably lies in the first months of 1913.\\
Many thanks in advance\\
Sincerely\\
Your obedient servant\\
L.E.J.~Brouwer.
\citep[Online Supplement, p.~1743, trl.~ours]{Dalen2011}
\end{quoting}
(`\dutch{Ondergeteekende zou het op hoogen prijs stellen,
indien hij den preciezen datum kon vernemen,
waarop Heft 2 van Band 142 van het
\german{Journal für die reine und angewandte Mathematik} (Crelle’s Journal)
ter Universiteits-bibliotheek is ontvangen.
De datum ligt waarschijnlijk in de eerste maanden van 1913.
Met beleefden dank bij voorbaat,
Hoogachtend,
Uw Dienstwillige L.E.J.~Brouwer.}')
The issue in question contained Brouwer's paper `\german{Über den
natürlichen Dimensionsbegriff}' \citep{Brouwer1913A}; the answer to his
question here is January 27, 1913 \citep[p.~358]{Dalen2008a}.} On the other
hand, we wish to note that neither in Brouwer's extensive notebooks
1904–1907, which show Brouwer to have been an omnivorous reader, nor
in his remaining correspondence, nor in the Dissertation have we found a
reference to the article by Pringsheim and Molk, nor to any other article
in the \emph{\french{Encylopédie}}, nor to any other of Molk's writings.

Also some younger French mathematicians were sensitive to the
issues later raised by Brouwer. For example, Lebesgue (1875–1941)
had stated, in a letter to Borel published in 1905:\footnoteA{\citet[p.~269]{Baire.Borel.Hadamard.Lebesgue1904}.}
\begin{quoting}
Although I strongly doubt that one will ever name a set that
is neither finite nor infinite, the impossibility of such a set
seems to me not to have been demonstrated.\footnoteA{\french{Bien
que je doute fort qu'on nomme jamais un ensemble qui ne soit
ni fini, ni infini, l'impossibilité d'un tel ensemble ne me
paraît pas démontré.}}
\end{quoting}
Note that one of the examples in Brouwer's paper of a principle that
has not been demonstrated is `every number is finite or infinite'. In
the intuitionistic setting, an example of a set that is neither finite
nor infinite was given in \citet[pp.~3–4]{Brouwer1924N}.

However, in spite of the early efforts by Kronecker, only with Brouwer
do we get a comprehensive development of mathematics excluding any
`unreliable' use of the principle of excluded middle.

\chapter{Direct influence}

In spite of its historical significance,
Brouwer's
paper has apparently had surprisingly little direct influence on others, apart
from sporadic references in Brouwer's own work
and that of Heyting.
The exception
is George François Cornelis Griss, as will now be explained.%
\footnoteA{The `Unreliability' paper may further have had an indirect influence already on Husserl.
In 1928 with his `\german{Intuitionistische Betrachtungen über den
Formalismus}' \citep{Brouwer1928A2},
Brouwer returned,
with explicit reference, to the themes of the earlier paper.
The relevance of Brouwer's 1928 paper for Husserl's
\emph{Formale und transzendentale Logik} \citep{Husserl1929}
is clear and,
several years ago,
was emphasized to one of us
by Byung-Hak Ha.
Afterwards Thomas Vongehr at the Husserl Archives in Louvain found
an entry in an old card catalogue that showed that Husserl had
owned an offprint of that paper.
Unfortunately,
the offprint itself was no longer to be found.
Brouwer and Husserl met
in
Amsterdam in April 1928
\citetext{\citealp[Online Supplement, p.~1515]{Dalen2011}; \citealp[vol.~5, p.~156]{Husserl1994b}},
and it is likely that Brouwer gave
the offprint to Husserl then,
or sent it in the aftermath.
However,
in spite of its topical closeness to some of the main
themes of
\emph{Formale und transzendentale Logik},
the latter contains no reference to Brouwer.
(We express our thanks to Ha and Vongehr.)}

There is a direct connection
between Brouwer's `Unreliability' and Griss' development, in a series of
papers published from 1944 to 1951, of a version of intuitionism without
negation.\footnoteA{\label{HeytingGriss}
\citet{Griss1944,Griss1946,Griss1950,Griss1951a,Griss1951b,Griss1951c}.
For more on
Griss and his work, see
\citet{Heyting1955} and \citet{Franchella1993}.}
Griss had first explained
his rationale to Brouwer directly, in a letter of April 19, 1941:
\begin{quoting}
Showing that something is not true, i.e.~showing the
incorrectness of a supposition is not an intuitively clear
act. For it is impossible to have an intuitively clear concept
of an assumption that later turns out to be even wrong. One
must maintain the demand that only building things up from the
foundations makes sense in intuitionistic mathematics.
\citep[p.~402, trl.~Van Dalen]{Dalen2011}\footnoteA{\dutch{Aantonen, dat iets niet waar is, d.w.z.~de onjuistheid van een veronderstelling aantonen, is niet een
intuïtief-duidelijke handelwijze. Van een veronderstelling,
die later zelfs blijkt fout te zijn, kan men namelijk onmogelijk
een intuïtief-duidelijke voorstelling hebben. Men moet de eis
handhaven, dat alleen het opbouwen vanaf de grondslagen in de
intuïtionistische wiskunde betekenis heeft.}
\citep[Online Supplement, p.~2142]{Dalen2011}}
\end{quoting}
This was repeated almost verbatim in Griss' publication of 1946:
\begin{quoting}
On philosophic grounds I think the use of the negation in
intuitionistic mathematics has to be rejected. Proving that
something is not right, i.e.~proving the incorrectness of
a supposition, is no intuitive method. For one cannot have a
clear conception of a supposition that eventually proves to
be a mistake. Only construction without the use of negation
has some sense in intuitionistic mathematics. \citep[p.~675]{Griss1946}
\end{quoting}
However, where in the article Griss prefers not to go further into the
philosophical issue and goes on to discuss mathematical consequences,
in the letter he first offers a justification for his basic idea. It
takes the form of a comment on Brouwer's `Unreliability':
\begin{quoting}
Although my ideas about the foundations of mathematics are not
completely identical to yours, the differences are unimportant
for what follows, so, for example, I can agree completely
with your considerations in the \emph{\dutch{Tijdschrift voor
Wijsbegeerte}}, 2nd volume, 1908. Let me just remark that the
concept of negation does not explicitly occur in the formulation
of the foundations of mathematics, but only in the examination
of the validity of the logical principles. You say there:
\begin{quoting}
The principle of contradiction is just as little in
dispute; the execution of the fitting of a system $a$
in a particular way into a system $b$, and finding that
this fitting turns out to be impossible are mutually
exclusive.\footnoteA{The translation of this passage that we will give below is a little different from this one by Van Dalen, but not substantially so. Note that another passage that Griss could have referred to is \citet[p.~127]{Brouwer1907}.}
\end{quoting}
What does impossibility of a `fitting in' mean here?

In the first place this can mean that one assumes the
possibility of fitting, and that this assumption leads to a
contradiction. This manner far exceeds the construction of
mathematical systems on the basis of the ur-intuition, and
as I remarked in the beginning, one cannot clearly obtain a
conception of it. If one still accepts it, then one takes in
principle a similar step, as when one accepts the principle of
the excluded third. An element of arbitrariness enters in our
idea about what is and what is not admissible in mathematics,
if one does not stick strictly to the requirement that one only
builds up mathematical systems from the foundations which are
given in the ur-intuition.

Another meaning which can be given to `finding that this fitting
of a system $a$ into a system $b$ turns out to be impossible'
might be this: that the system $a$ demonstrably differs (in that
case this concept has to be defined) from every system that can be
fitted into $b$. One asks for example whether $e$ is an algebraic
number and one finds that $e$ is positively transcendent so $e$
demonstrably differs from each algebraic number. If need be,
one can even answer the question whether $e$ is algebraic by:
$e$ is not algebraic, but then we have assigned a new meaning
to the word `not'.\footnoteA{%\vspace{-\baselineskip}
%\begin{quoting}\footnotesize
\dutch{Hoewel mijn ideeën over de grondslagen van
wiskunde niet volkomen gelijk zijn aan de Uwe, zijn de
verschillen voor het volgende niet van belang, zodat
ik bijv.~geheel kan aansluiten bij Uw beschouwingen
in het \emph{Tijdschrift voor Wijsbegeerte}, 2de jaargang,
1908. Alleen merk ik op, dat het begrip negatie bij het
formuleren van de grondslagen der wiskunde niet expliciet
optreedt, maar pas bij het onderzoek naar de geldigheid
der logische principes. U zegt daar:
\begin{quoting}[leftmargin=0em,rightmargin=0em]
Evenmin is aanvechtbaar het principe van
contradictie: het volvoeren van de inpassing
van een systeem a op een bepaalde wijze in een
systeem b, en het stuiten op de onmogelijkheid
van die inpassing sluiten elkaar uit.
\end{quoting}
}
\noindent\dutch{Wat betekent hierin de onmogelijkheid
van een inpassing?}

\hspace{\footparindent}\dutch{Ten eerste kan dit betekenen, dat men
van de mogelijkheid van een inpassing uitgaat en
die veronderstelling ad absurdum voert. Deze wijze
van doen gaat echter ver uit boven het opbouwen van
wiskundige systemen op grond van de oerintuïtie, en
men kan er zich, zoals ik reeds in het begin opmerkte,
geen duidelijke voorstelling van maken. Aanvaardt men
ze toch, dan doet men in principe een dergelijke stap,
als wanneer men het beginsel van het uitgesloten derde
zou aanvaarden. Er komt een element van willekeur in
onze opvatting over wat al of niet toelaatbaar is in de
wiskunde, als men niet streng vasthoudt aan de eis alleen
wiskundige systemen op te bouwen vanaf de grondslagen,
die in de oerintuïtie gegeven zijn.}

\hspace{\footparindent}\dutch{Een andere betekenis, die aan de uitdrukking
`stuiten op de onmogelijkheid van een inpassing van een
systeem a in een systeem b' gehecht kan worden, zou deze
kunnen zijn, dat het systeem a aanwijsbaar verschilt (dit
begrip moet dan gedefinieerd zijn) van ieder systeem,
dat in het systeem b kan worden ingepast. Men vraagt
bijvoorbeeld, of $e$ een algebraïsch getal is en vindt,
dat $e$ positief transcendent is, zodat $e$ aanwijsbaar
van ieder algebraïsch getal verschilt. Desnoods kan
men zelfs op de vraag, of $e$ algebraïsch is, antwoorden:
$e$ is niet algebraïsch, maar dan hebben we aan het woord
`niet' een nieuwe betekenis toegekend.}%\vspace{-\baselineskip}
%\end{quoting}
}
% end footnote
\end{quoting}
Brouwer's paper `Essentially negative properties' \citep{Brouwer1948A} was written
in response to Griss. In his letter of 1941, Griss had remarked that
`no real number $a$ is known about which it has been proved that it
cannot possibly be equal to 0 ($a \neq 0$), while at the same time it
has not been proven that the number differs positively from 0 
$(a\mathrel{\#}0)$'. Brouwer in his paper constructed a real number $a$ with just that
property; but he did not provide an accompanying philosophical account
as an alternative to Griss' view.

An occasion for Brouwer, Griss and others to debate these
matters in public would have been a meeting planned by S.I.~Dockx.\footnoteA{Stanislas Isnard Dockx, Antwerp 1901–Brussels 1985.}
% , \url{http://www.lesacademies.org/iist/biographie.html
A letter of Beth
to Dockx of July 8, 1949, suggests that also Freudenthal, Heyting,
and Van Dantzig were invited, but at the same time makes it clear that
Brouwer declined because he did not want to participate in an event
with Freudenthal \citep[Online Supplement, p.~2446]{Dalen2011}.\footnoteA{For an account of
Brouwer's by then long-standing conflict with Freudenthal, see \citet[pp.~721–728, 753–757, and 794–799]{Dalen2005}.}
To the best of our knowledge, the meeting never took place.

Heyting published a reaction in 1955, `G.F.C.~Griss and his negationless
intuitionistic mathematics'. While Heyting noted that `unrealized suppositions' are implicit
in all general statements, so that banishing such suppositions would
reduce mathematics to an `utterly unimportant and uninteresting subject'
\citep[p.~95]{Heyting1955}, he did not provide a detailed confrontation with the arguments
of Griss. It can be argued that Brouwer's dissertation in
effect contains an answer to Griss' objection:
according to
\citet{Atten2009b},
the view expressed on the hypothetical judgement at the beginning of chapter 3 of Brouwer's thesis
\citep[pp.~125–127]{Brouwer1907}
is that
logical reasoning does not
operate on constructions, let alone hypothetical ones, but on conditions
on constructions. The difference is that these conditions, whether
fulfillable or not, can themselves be represented as actual objects.

% needed for pages environment;
% unlike article.cls, thpl.cls does not have quotation

\let\quotation\quote
\let\endquotation\endquote

\begin{pages}

\begin{Leftside}

\selectlanguage{dutch}

\beginnumbering

\autopar

% title
\pstart
\title{De onbetrouwbaarheid der logische principes}
\author{door\\ L.E.J.~Brouwer}
\date{}
\emptythanks
\maketitle
%{\noindent\huge\center De onbetrouwbaarheid der logische principes\par}
%{\noindent\Large\center door\\ L.E.J.~Brouwer\par}
%\pend[\vspace{2\baselineskip}]
%\vspace{2\baselineskip}
%\pend[\vspace{2\baselineskip}]
%~\\
%~\\
%\pstart
%\noindent
{}[Noot Brouwer voorafgaand aan de herdruk van 1919]
\emph{Dit opstel zou ook thans nog
in denzelfden vorm geschreven kunnen zijn.
Medestanders hebben de er verdedigde opvattingen
nog weinig gevonden.}
\pend[\vspace{1\baselineskip}]

%2
1.
De
\emph{wetenschap}
beschouwt herhaling in den tijd van als
onderling gelijk stelbare volgreeksen van qualitatieve
verscheidenheid in den tijd.
Dit vereenzamen der idee tot waarneembaarheid,
en als zoodanig tot herhaalbaarheid,
verschijnt na religielooze scheiding%
\footnoteB{een vermogen,
voortgekomen uit de oerzonde van vrees of begeerte,
maar wederkeerend,
ook zonder levende vrees of begeerte.
vgl.~L.E.J.~Brouwer. Leven, Kunst en Mystiek. pag.~13–23.}
tusschen subject en
tot \emph{iets anders} geworden onbereikte bereikbaarheid.
De drang tot bereiking dezer bereikbaarheden wordt in het
intellect volgens een wiskundig systeem van gestelde
stelbaarheden,
geboren uit abstractie van herhaling en herhaalbaarheden,
gestuurd langs onmiddellijke bereiktheden.

%3
Alles wat verschijnen kan als onbereikte bereikbaarheid,
laat zich in systemen van gesteldheden intelligeeren,
zoo ook religie;
maar dan is de religieuze \emph{wetenschap} religieloos:
gewetensussend,
of ijdel spel,
of slechts van doelnajagende beteekenis.%
\footnoteB{t.a.p.~pag.~27.}

%4
En,
als alle religieloosheid,
heeft wetenschap nòch religieuze betrouwbaarheid,
nòch betrouwbaarheid in zich.
\pagenumber{153}
Allerminst kan een wiskundig systeem van gesteldheden,
los van de waarnemingen,
die het
intelligeerde,
onbepaald vervolgd,
betrouwbaar blijven in het richten langs die waarnemingen.

%5
Zoodat onafhankelijk van de waarneming volvoerde logische
redeneeringen,
die immers beteekenen wiskundige transformaties in het
intelligeerende wiskundig systeem,
uit wetenschappelijk aanvaarde praemissen onaannemelijke
conclusies kunnen afleiden.%
\footnoteB{t.a.p.~pag.~20, 21.}

%6
De klassieke opvatting,
die in de ervaringsgeometrie uit aanvaarde praemissen door
volgens de logische principes gevoerde redeneeringen slechts
onaanvechtbare conclusies zag afleiden,
induceerde de logische redeneeringen als methode van opbouw
der wetenschap en de logische principes als menschelijke
vermogens tot opbouw van wetenschap.

%7
Maar de geometrische redeneeringen gelden slechts voor een
onafhankelijk van eenige ervaring in het intellect
opbouwbaar wiskundig systeem,
en dat een zoo populaire groep van waarnemingen als de
geometrie het bedoelde wiskundig systeem zoo blijvend
verdraagt,
verdient,
als alle proefhoudende natuurwetenschap,
met wantrouwen te worden aangezien.

%8
Het inzicht van de wetenschappelijke onbetrouwbaarheid der
logische redeneeringen maakt,
dat de conclusiën van Aristoteles omtrent de constitutie der
natuur zonder practische verifieering niet overtuigen;
dat de waarheid,
die bij Spinoza opengaat,
geheel onafhankelijk wordt gevoeld van zijn logische
systematiek;
dat men niet gehinderd wordt door de antinomieën van Kant,
en evenmin door het ontbreken van in al haar consequenties
door te voeren physische hypothesen.

%9
\pstart
\pagenumber{154} Bovendien zijn bij de betoogen betreffende op wiskundige
systemen gespannen ervaringswerkelijkheden de logische
principes niet het richtende,
maar in de begeleidende taal achteraf opgemerkte
regelmatigheid,
en zoo men los van wiskundige systemen spreekt volgens die
regelmatigheid,
is er altijd gevaar voor paradoxen als die van Epimenides.
\pend[\vspace{\baselineskip}]

%10
\pstart
2. In religieuze waarheid,
in \emph{wijsheid},
die de splitsing opheft in subject en iets anders,
is geen wiskundig intelligeeren,
daar de verschijning van den tijd niet langer wordt
aanvaard,
nog minder dus betrouwbaarheid van logica.
Integendeel,
de taal der inkeerende wijsheid verschijnt ordeloos,
onlogisch,
omdat ze nooit kan voeren langs in het leven gedrukte
systemen van gesteldheden,
slechts hun breking kan begeleiden,
en zoo misschien de wijsheid,
die die breking doet,
kan laten opengaan.%
\footnoteB{t.a.p.~pag.~47 vlgg.,
65 vlgg.}
\pend[\vspace{\baselineskip}]

%11
3.
Blijft de vraag,
of dan althans de logische principes vaststaan voor van
levensinhoud vrije \emph{wiskundige systemen},
voor systemen opgetrokken uit de gestelde abstractie van
herhaling en herhaalbaarheid,
uit de gestelde inhoudslooze tijdsintuïtie,
uit de oer-intuïtie der wiskunde.%
\footnoteB{vgl.~L.E.J.~Brouwer.
Over de Grondslagen der Wiskunde. pag.~8, 81, 98,
179.}
Door alle tijden is in wiskunde met vertrouwen logisch
geredeneerd;
nooit aarzelde men,
door logica uit postulaten getrokken conclusies te
aanvaarden,
waar de postulaten gelden.
In dezen tijd zijn echter paradoxen geconstrueerd,
die wiskundige
paradoxen schijnen%
\footnoteB{Burali-Forti. (\italian{Rendiconti del circolo\
Matematico di Palermo}. 1897. p.~164).\\
Zermelo. (\german{Mathematische Annalen} 59). Koenig. (ibid.~61).\\
Richard. (\french{Revue générale des Sciences}. 1905).\\
Russell. (The Principles of Mathematics. Part I. Chap.~X).\\
Voor pogingen tot oplossing dezer paradoxen vgl.,
behalve de opstellers zelf:
Poincaré. (\french{Revue de Métaphysique et de Morale}. 1905
no.~6, 1906 no.~1, 3). Mollerup. (\german{Mathematische Annalen}
64). Schoenflies. (\german{Bericht über die Mengenlehre}. \textsc{ii}\@. Kap.~\textsc{i}.~\S 7).},
en wantrouwen wekken tegen het
\pagenumber{155}
vrije gebruik van logica in
wiskunde,
zoodat enkele wiskundigen hun vooronderstelling van logica
in wiskunde loslaten,
en logica en wiskunde tezamen trachten op te
bouwen%
\footnoteB{in het bijzonder Hilbert in
\german{Verhandlungen des internationalen Mathematiker-Congresses
in Heidelberg 1904}. p.~174.},
in aansluiting aan de door Peano gegrondveste school der
\emph{logistiek}.
Aangetoond kan echter worden%
\footnoteB{Grondslagen der Wiskunde. \textsc{iii}.},
dat deze paradoxen voortkomen uit dezelfde dwaling als die
van Epimenides,
dat ze namelijk ontstaan,
waar regelmatigheid in de taal,
die wiskunde begeleidt,
wordt uitgebreid over een taal van wiskundige woorden,
die geen wiskunde begeleidt;
dat verder de logsitiek eveneens zich bezighoudt met de
wiskundige taal in plaats van met de wiskunde zelf,
dus de wiskunde zelf niet verheldert;
dat ten slotte alle paradoxen verdwijnen,
als men zich beperkt,
slechts te spreken over expliciet uit de oer-intuïtie
opbouwbare systemen,
m.a.w.~in plaats van logica door wiskunde,
wiskunde door logica laat vooronderstellen.

%12
Zoo blijft nu alleen nog de meer gespecialiseerde vraag:
“Kan men bij zuiver wiskundige constructies en
transformaties de voorstelling van het opgetrokken wiskundig
systeem tijdelijk verwaarloozen,
en zich
bewegen in het accompagneerend taalgebouw,
geleid door de principes van \emph{syllogisme},
van \emph{contradictie} en van \emph{tertium exclusum},
in vertrouwen dat door tijdelijke oproeping van de
voorstelling der beredeneerde wiskundige constructies
telkens elk deel van het betoog zou kunnen worden
gewettigd?”

%13
Hier zal blijken,
dat dit vertrouwen voor de beide eerste principes wèl,
voor het laatste niet gegrond is.

%14
Het \emph{syllogisme} vooreerst leest in de inpassing van een
\pagenumber{156}
systeem
$b$
in een systeem
$c$
en de daarmee samengaande
inpassing van een systeem
$a$
in het systeem $b$
een directe
inpassing van het systeem
$a$
in het systeem
$c$,
wat niet anders is dan een tautologie.

%15
Evenmin is aanvechtbaar het principe van \emph{contradictie}:
het volvoeren van de inpassing van een systeem
$a$
op bepaalde
wijze in een systeem
$b$,
en het stuiten op de onmogelijkheid van die inpassing
sluiten elkander uit.

%16
Nu het principium \emph{tertii exclusi}:
dit eischt,
dat iedere onderstelling òf juist òf onjuist is,
wiskundig:
dat van iedere onderstelde inpassing van systemen op
bepaalde wijze in elkaar hetzij de beëindiging,
hetzij de stuiting op onmogelijkheid kan worden
geconstrueerd.
De vraag naar de geldigheid van het principium tertii
exclusi is dus aequivalent met de vraag naar de
\emph{mogelijkheid van onoplosbare wiskundige problemen}.
Voor de wel eens uitgesproken%
\footnoteB{vgl.~Hilbert.
\german{Mathematische Probleme.
Göttinger Nachrichten.
1900}.
Ook Schoenflies
(l.c.)
wil onvoorwaardelijk de methode van het indirecte bewijs
handhaven,
die hij ten onrechte uitsluitend van het principium
contradictionis afhankelijk acht.}
overtuiging,
dat onoplosbare wiskundige
problemen niet bestaan,
is geen aanwijzing van een bewijs aanwezig.

%17
Zoolang alleen bepaalde eindige discrete systemen gesteld
worden,
is het onderzoek naar de mogelijkheid of onmogelijkheid
eener inpassing steeds beëindigbaar en voerend tot antwoord,
is dus het principium tertii exclusi een betrouwbaar
redeneerprincipe.%
\footnoteB{Dit onderzoek kan zelfs steeds
door een machine worden uitgevoerd,
of door een gedresseerd dier,
vereischt niet de oer-intuïtie der wiskunde,
levend in een menschelijk intellect.
Maar tegenover vragen betreffende oneindige verzamelingen
wordt die oer-intuïtie telkens weer onmisbaar;
door dit voorbij te zien,
zijn Peano en Russell,
Cantor en Bernstein slechts tot dwalingen gekomen.}

%18
Dat ook oneindige systemen ten opzichte van zoovele
eigenschappen eindig worden beheerscht,
geschiedt door overzien van de aftelbaar oneindige reeks der
geheele ge\pagenumber{157}tallen met
\emph{volledige inductie}%
\footnoteB{Poincaré is
misschien de eenige,
die in de volledige inductie
`\french{le raisonnement mathématique
par excellence}' heeft herkend.
Vgl.~La Science et l'Hypothese.
Chap.~I.},
namelijk door opmerken van eigenschappen,
d.w.z.~inpassingen,
die voor een
\emph{willekeurig geheel getal}
gelden,
in het bijzonder ook van contradicties,
dat zijn onmogelijke inpassingen,
die voor een
\emph{willekeurig geheel getal}
gelden.
Dat echter uit de in een vraag gestelde systemen een is af
te leiden,
dat door een invariant over een aftelbaar oneindige reeks de
vraag volledig induceerend leest,
en zoo oplost,
blijkt eerst a posteriori,
als toevallig de constructie van zulk een systeem gelukt is.
Want het geheel der uit de vraagstelling te ontwikkelen
systemen is
\emph{aftelbaar onaf}%
\footnoteB{vgl.~\dutch{Grondslagen der Wiskunde}. p.~148.},
dus niet a priori methodisch te
onderzoeken ten opzichte van de aanwezigheid of afwezigheid
van een de vraag beslissend systeem.
En het is niet uitgesloten,
dat een even gelukkige greep,
als zoo dikwijls de beslissing bracht,
eens het aftelbaar onaffe systeem der mogelijke
ontwikkelingen tot een onoplosbaarheid zou overzien.

%19
Zoodat in oneindige systemen het principium tertii exclusi
vooralsnog niet betrouwbaar is.
Toch zal men bij ongerechtvaardigde toepassing nooit kunnen
stuiten op een contradictie en zoo de ongegrondheid van zijn
redeneeringen ontdekken.
Immers daartoe zouden de volvoering en de contradictoriteit
van een inpassing beide tegelijk contradictoor moeten kunnen
zijn,
wat het principium contradictionis niet toelaat.

%20
Een sprekend voorbeeld levert de volgende onbewezen
stelling,
die op grond van het principium tertii exclusi in de
gangbare theorie der transfinite getallen algemeen vertrouwd
en gebruikt wordt,
dat n.l.~elk getal is òf eindig òf oneindig,
m.a.w.~dat voor elk getal $\gamma$ kan worden geconstrueerd:

%21
\pagenumber{158}
hetzij een afbeelding van $\gamma$ geheel op de rij der
geheele getallen zóó,
dat daarbij een getal $\alpha$ uit die
rij
\emph{het laatste}
is (de getallen $\alpha+1, \alpha+2,
\alpha+3, \dots$ vrij blijven),

%22
hetzij een afbeelding van $\gamma$ geheel of gedeeltelijk op
de rij der geheele getallen in haar geheel.%
\footnoteB{De
eventueele onjuistheid dezer stelling zal weer nooit in een
contradictie kunnen blijken; immers de contradictoriteit van
de constructie der vrij blijvende rij
$\alpha+1, \alpha+2, \alpha+3 \dots$,
en die van haar contradictoriteit kunnen
nooit tezamen optreden.}

%23
Zoolang deze stelling onbewezen is,
moet men voor onzeker houden,
of vragen als:

%24
“\emph{Is bij de decimale ontwikkeling van $\pi$ een cijfer,
dat duurzaam veelvuldiger optreedt,
dan alle andere?}”

%25
“\emph{Komen bij de decimale ontwikkeling van $\pi$ oneindig
veel paren van gelijke opeenvolgende cijfers voor?}”

%26
een oplossing bezitten.

%27
En evenzoo onzeker blijft,
of de algemeenere wiskundige vraag:

%28
“\emph{Is in de wiskunde het principium tertii exclusi onbepaald
geldig?}”

%29
\pstart
een oplossing bezit.%
\footnoteB{Men behoort dus in wiskunde
de gewoonlijk als
\emph{bewezen}
geldende stellingen te
onderscheiden in
\emph{juiste}
en
\emph{niet-contradictore}.
Tot de eerste behooren de algebraïsche en analytische
gelijkheden,
en de geometrische snijpuntsstellingen;
ook,
dat een puntverzameling geen andere machtigheid bezitten
kan,
dan de (\dutch{Grondslagen}, pag.~149) genoemde.
Tot de laatste,
dat een puntverzameling zeker een dier machtigheden bezit;
ook,
dat een afgesloten puntverzameling zich laat splitsen in een
perfekte en een aftelbare.}
\pend[\vspace{\baselineskip}]

%30
Samenvattende:

%31
In wijsheid is geen logica.

%32
In wetenschap is logica vaak,
maar niet duurzaam doeltreffend.

%33
In wiskunde is niet zeker,
of alle logica geoorloofd is,
en is niet zeker of is uit te maken,
of alle logica geoorloofd is.

\endnumbering

\end{Leftside}

\begin{Rightside}

\selectlanguage{english}

\beginnumbering

\autopar

% title
\pstart
\title{The unreliability of the logical principles}
\author{by\\ L.E.J.~Brouwer}
\date{}
\emptythanks
\maketitle
%{\noindent\huge\center \par}
%{\noindent\Large\center  \par}
%\pend[\vspace{2\baselineskip}]
%~\\
%~\\
%\pstart
%\noindent
{}[Brouwer's note preceeding the 1919 reprint]
\emph{This essay could also today still
be written in the same form. The opinions defended in it have, as yet,
found few supporters.}
\pend[\vspace{1\baselineskip}]

%2
1. \emph{Science} considers repetition in time of interidentifiable
succession-sequences of qualitative differentiation through time. This
isolating of the idea into an observable, and as such a repeatable,
emerges after areligious\footnoteD{NB not `antireligious'.}
separation\footnoteC{A capacity, rooted in the original sin of fear
or desire, but reappearing also without living fear or desire. Cf.~\citet[pp.~13–23]{Brouwer1905A}.} between
the subject and an unreached reachable that has become \emph{something
separate}. In the intellect, the urge to reach these reachables is
conducted along things immediately reached, according to a mathematical
system of posited positables, born out of abstraction of repetition
and repeatables.

%3
Everything that can emerge as unreached reachable 
lets itself be intelligized
in systems of posits,\footnoteD{`Posit' not in Quine's sense,
but rather like Kant's
`\german{Setzung}'.} thus also religion; but then religious
\emph{science} is areligious:
conscience-numbing, or idle play, or of merely goal-chasing
significance.\footnoteC{\citet[p.~27]{Brouwer1905A}.}

%4
And science, as everything areligious, possesses \emph{neither} religious
reliability, \emph{nor} reliability in itself.
Least of all can a mathematical system of posits, separated from the
observations it made intelligible, when continued indefinitely, remain
reliable when directing along those observations.

%5
Consequently logical argumentations, which, after all, consist in
mathematical transformations in the mathematical system that makes
[the observations] intelligible, may derive unlikely conclusions from
scientifically accepted premises, when carried out independently of
observation.\footnoteC{\citet[pp.~20, 21]{Brouwer1905A}.}

%6
The classical conception, which in experiential geometry witnessed
reasonings \textemdash\ from accepted premises, carried out according to
logical principles \textemdash\ derive only incontrovertible conclusions,
induced the logical reasonings as method for building science, and the
logical principles as human capacities for building science.

%7
But the geometrical reasonings are valid only for a mathematical system
that can be built in the intellect independently of any experience,
and that such a popular group of observations as geometry corroborates
the mathematical system in question so enduringly, deserves, like
all experimental natural science, to be regarded with distrust.

%8
The
insight of the scientific unreliability of the logical reasonings has as
consequence that Aristotle's conclusions on the constitution of nature are
unconvincing without verification in practice; that the truth unveiled in
Spinoza is experienced wholly independently of his logical architectonic;
that one is not hindered by the antinomies of Kant, nor by the lack of
physical hypotheses that can be carried through in all their consequences.

%9
\pstart
Moreover, regarding discourse concerning experiential realities that
have been cast on mathematical systems, the logical principles are not
directive, but regularities that have afterward been noticed in the
accompanying language, and if one speaks according to these regularities
with no link to mathematical systems, there is always the danger of
paradoxes such as that of Epimenides.
\pend[\vspace{\baselineskip}]

%10
\pstart
2. In religous truth, in \emph{wisdom}, which suspends the splitting into
subject and something separate, there is no mathematical intellection, as
the appearance of time is no longer accepted, even less thus the reliability of logic. On the contrary,
the language of inward-turning wisdom appears without order, illogical,
because it can never carry along systems of posits pressed upon life,
but can only accompany their breakdown, and thus perhaps unveil
the wisdom that effects the break.\footnoteC{\citet[p.~47ff,65ff]{Brouwer1905A}.}
\pend[\vspace{\baselineskip}]

%11
3. The question remains whether then the logical principles hold
at least for \emph{mathematical systems} that are free of living
content,\footnoteD{For the expression `living content', compare:
`… the 
\emph{existence} of that mathematical reasoning system does
not entail that it 
\emph{lives},
in other words that it accompanies a chain of thoughts …'
\citet[p.~138n, emphasis Brouwer, trl.~ours]{Brouwer1907}
(\dutch{… volgt uit het
\emph{bestaan}
van dat wiskundig redeneersysteem
nog niet,
dat dat taalsysteem
\emph{leeft},
m.a.w.~een aaneenschakeling van gedachten
begeleidt\dots}).} for systems erected from posited
abstraction of repetition and repeatability, from the posited contentless
intuition of time, from the [Ur-]intuition of mathematics.\footnoteC{Cf.~\citet[pp.~8, 81, 98,179.]{Brouwer1907}.} Through all ages, in mathematics one has reasoned logically
with confidence; never did one hesitate to accept conclusions drawn
from postulates by logic, where the postulates hold. In this time,
however, paradoxes have been constructed that seem to be mathematical
paradoxes,\footnoteC{\citet{Burali-Forti1897},
\citet{Zermelo1904},
\citet{Koenig1905},
\citet{Richard1905},
\citet[Part I. Chap.~X]{Russell1903}.
For attempts at solving these paradoxes see, besides the proposers
themselves:
\citet{Poincare19051906},
\citet{Mollerup1907},
\citet[Kap.~I.~\S~7]{Schoenflies1908}.
[Note that the preface to that work is dated `\german{im Oktober 1907}'.]}
and that arouse distrust against the free use of logic in
mathematics, so that some mathematicians let go of their presupposition
of logic in mathematics, and try to build up logic together with
mathematics,\footnoteC{In particular Hilbert in \citet{Hilbert1905}.} following the school of \emph{logistics} founded by Peano. It can
be shown,%
\footnoteC{\citet[ch.~\textsc{iii}]{Brouwer1907}.}
however, that these paradoxes result from the same error as
that of Epimenides, namely, that they arise where regularities in the
language that accompanies mathematics are extended over a language of
mathematical words that does not accompany mathematics; that, further,
logistics too is concerned with the mathematical language instead of
with mathematics itself, thus does not clarify mathematics itself; that,
finally, all paradoxes disappear, when one restricts oneself to speaking
only of systems that explicitly can be built out of the Ur-intuition,
in other words, when instead of letting mathematics presuppose logic,
one lets logic presuppose mathematics.

%12
Thus, now only the more specific question still remains: `Can one,
in the case of purely mathematical constructions and transformations,
temporarily neglect the presentation of the mathematical system that has
been erected, and move in the accompanying linguistic building, guided by
the principles of \emph{the syllogism}, of \emph{contradiction}, and of
\emph{tertium exclusum}, always confident that, by momentary evocation
of the presentation of the mathematical constructions suggested by this
reasoning, each part of the discourse could be justified?'

%13
Here it will turn out that this confidence is well-founded for each
of the first two principles, but not for the last.

%14
To begin with, the
\emph{syllogism} reads the fitting of a system \emph{b} into a system
\emph{c} and the concommitant fitting of a system \emph{a} into the
system \emph{b} as a direct fitting of the system \emph{a} into the
system \emph{c}, which is nothing but a tautology.\footnoteD{A tautology
in the sense of, for example, Kant, not Wittgenstein.}

%15
Nor can the principle of \emph{contradiction} be assailed: completing the
fitting of a system \emph{a} in a certain way into a system \emph{b}, and
being blocked by the impossibility of that fitting, exclude one another.

%16
Now the principium \emph{tertii exclusi}: this demands
that every
supposition is either correct or incorrect, mathematically: that of
every supposed fitting in a certain way of systems in one another,
either the termination or the blockage by impossibility, can be
constructed. The question of the validity of the principium tertii exclusi
is thus equivalent to the question concerning the \emph{possibility of
unsolvable mathematical problems}. For the already proclaimed conviction
that unsolvable mathematical problems do not exist, no indication of a
demonstration is present.\footnoteC{Cf.~\citet{Hilbert1900}. Also Schoenflies \citeyearpar{Schoenflies1908} wants to
uphold the method of indirect proof unconditionally, which he mistakenly
considers to depend only on the principium contradictionis.}

%17
As long as only certain finite discrete systems are posited, the
investigation into the possibility or impossibility of a fitting can
always be terminated and leads to an answer, whence the principium tertii
exclusi is a reliable principle of reasoning.\footnoteC{This investigation
itself can always be done by a machine or by a trained animal, not
requiring the intuition of mathematics living in a human intellect. But
in face of questions involving infinite sets, that intuition becomes,
again and again, indispensable; by overlooking this, Peano and Russell,
Cantor and Bernstein have only arrived at errors. [Brouwer gives an
exposition of these errors in the chapter 3
of his dissertation \citealt{Brouwer1907}, `\dutch{Wiskunde en Logica}'.]}

%18
That also infinite systems, with respect to so many properties,
are controlled by finite means, is achieved by surveying the
denumerably infinite sequence of the whole numbers by \emph{complete
induction},\footnoteC{Poincaré is perhaps the only one who has recognized
mathematical induction as `\french{le raisonnement mathématique par
excellence}'. See \citet[Chap.~I]{Poincare1902}.}
namely by observing properties, that is, fittings, that hold for an
\emph{arbitrary whole number}, and in particular also contradictions,
that is, impossible fittings, that hold for an \emph{arbitrary whole
number}. However, that from the systems posited in a question, one can
be derived that reads the question by means of a complete induction,
on the basis of an invariant in a denumerably infinite sequence,
and thereby solves it, is found only a posteriori, when accidentally
the construction of such a system has succeeded. For the whole of the
systems that can be developed from the question posed is \emph{denumerably
unfinished},\footnoteC{Cf.~\citet[p.~148]{Brouwer1907}.}
whence cannot be a priori investigated methodically regarding the presence
or absence of a system that decides the question. And it is not excluded,
that by a draw as lucky as the ones that have so often led to a decision,
we will one day see from the denumerably infinite system of possible
developments that it is unsolvable.

%19
So that in infinite systems the principium tertii exclusi is as yet not
reliable. Still, one can never, in unjustified application, be blocked
by a contradiction and thereby discover the groundlessness of one's
reasonings. After all, to that end it would have to be possible for
the execution and contradictoriness of a fitting to be simultaneously
contradictory, which the principium contradictionis does not allow.

%20
A striking example is provided by the following undemonstrated
proposition, which on the ground of the principium tertii exclusi is
generally trusted and used in the current theory of transfinite numbers,
namely that every number is either finite or infinite, in other words,
that for every number $\gamma$ one can construct:

%21
either a mapping of all of $\gamma$ to the sequence of the whole numbers,
such that a number $\alpha$ in that sequence is \emph{the last one}
(the numbers $\alpha+1, \alpha+2, \alpha+3, \dots$ remain free),

%22
or a mapping of all or part of $\gamma$ to the sequence of the
natural numbers in its entirety.\footnoteC{A latent incorrectness of
this proposition also shall never become clear from a contradiction:
after all, the contradictoriness of the construction of the sequence
$\alpha+1, \alpha+2, \alpha+3, \dots$ which remains free and that of
its contradictoriness can never occur together.}

%23
As long as this proposition is undemonstrated, it must be held uncertain
whether questions such as:

%24
`\emph{Is there in the decimal expansion of $\pi$ a digit that occurs
enduringly more often than all others?}'

%25
`\emph{Do there occur in the decimal expansion of $\pi$ infinitely many
pairs of equal consecutive digits?}'

%26
have a solution.

%27
And likewise, it remains uncertain whether the more general mathematical
question:

%28
`\emph{Is in mathematics the principium tertii exclusi unconditionally valid?}'

%29
\pstart
has a solution.%
\footnoteC{One should therefore in mathematics distinguish
the propositions that are usually taken to have been \emph{demonstrated}
into \emph{correct} and \emph{non-contradictory} ones. To the former
belong the algebraic and analytic equalities, and the geometrical
incidence theorems; also, that a point set can have no other cardinality
than those mentioned in \citet[p.~149]{Brouwer1907}. To the latter,
that a point set does indeed have one of those cardinalities; also,
that a closed point set can be split into a perfect and a denumerable
one.}
\pend[\vspace{1\baselineskip}]

%30
Summarizing:

%31
In wisdom is no logic.

%32
In science, logic is often, but not enduringly efficacious.

%33
In mathematics it is not certain whether all logic is permissible,
and it is not certain whether it can be decided whether all logic is
permissible.

\endnumbering

\end{Rightside}

\end{pages}
\Pages

% whatever comes after (Acknowledgement, references)

\chapter*{Acknowledgements}

We thank the director of the Brouwer Archive,
Dirk van Dalen,
for permission to publish a translation
of Brouwer's article.
A French version
of our introduction,
translated from the original English by
Mark van Atten,
Vanessa van Atten,
and Michel Bourdeau,
together with a French translation of Brouwer's
paper by Mark van Atten and Michel Bourdeau,
was published
first,
as
Mark van Atten,
Göran Sundholm,
Michel Bourdeau et Vanessa van Atten
(2014)
\french{‘~«~Que les principes de la logique ne sont pas fiables.~»
Nouvelle traduction française annotée et commentée de
l’article de 1908 de L.E.J.~Brouwer’,
\emph{Revue d’histoire des sciences}} 67(2),
257–281.
We thank its publisher,
Armand Colin (Paris),
for permission to publish
here what is in effect a
(somewhat modified)
translation of that introduction.
The work with Michel Bourdeau and Vanessa van Atten on
the French version also led to some revisions in our
original English,
for which we are thankful to them.
We acknowledge generous financial support from the \textsc{ihpst}
(Paris),
assured by
Michel Bourdeau and Jean Fichot,
that twice enabled Sundholm to visit
Paris and made our joint work possible.

We are grateful to
Dirk van Dalen,
Per Martin-Löf,
and
Robert Tragesser
for comments on an earlier version of our translation,
and to Arthur Schipper
for suggesting various improvements to the manuscript.

\bibliographystyle{plainnat}

% see http://tex.stackexchange.com/questions/40747/bibtex-handling-of-the-dutch-van-name-prefix-with-natbib

\DeclareRobustCommand{\VAN}[3]{#3}

\bibliography{Brouwer1908NL-EN}

\end{document}